\documentclass[11pt]{amsart}
\usepackage[active]{srcltx}
\usepackage{pdfsync}
\usepackage{amsfonts,amssymb}
\usepackage{amsmath}
\usepackage{tabularx} \usepackage{arydshln}
\usepackage{wasysym}

\usepackage[kerning=true]{microtype}

\usepackage{xcolor}
\usepackage{a4wide}

\newtheorem{theorem}{Theorem}
\newtheorem{lemma}{Lemma}

\newtheorem{remark}{Remark}

\newcommand{\CC}{\mathbb{C}}

\newcommand{\RR}{\mathbb{R}}
\newcommand{\esf}{\mathbb{S}}
\newcommand{\esfera}{\mathbb{S}^{2n+1}}

\newcommand{\Ric}{\mathrm{Ric}}

\newcommand{\SU}{\mathrm{SU}}

\renewcommand{\hom}{\mathrm{Hom}}
\renewcommand{\Re}{\mathrm{Re}}
\def\dim{\mathop{\hbox{\rm dim}}}

\newcommand{\sof}{\mathfrak{so}}

\newcommand{\mm}{\mathfrak{m}}
\newcommand{\hh}{\mathfrak{h}}
\newcommand{\g}{\mathfrak{g}}

\begin{document}

\title[ Einstein connections  on Berger spheres]{  Einstein with skew-torsion connections \\on Berger spheres}
%\author{Cristina~Draper, Antonio~Garvín, Francisco~J.~Palomo}

\author[C.~Draper]{Cristina Draper Fontanals${}^*$}
\address{Departamento de Matem\'{a}tica Aplicada,   Universidad de M\'{a}laga,
  29071 M\'{a}laga, Spain}
\email{cdf@uma.es}
\thanks{${}^*$ Supported by the Spanish Ministerio de Econom\'{\i}a y Competitividad---Fondo Europeo de
Desarrollo Regional (FEDER) MTM2016-76327-C3-1-P, and by the Junta de Andaluc\'{\i}a grants FQM-336 and FQM-7156, with FEDER funds}

\author[A.~Garv\'{\i}n]{Antonio Garv\'{\i}n${}^\circ$}
%\curraddr{}
\email{garvin@uma.es}
\thanks{${}^\circ$ Supported by the
Spanish MEC grant MTM2013-41768-P and by the Junta de Andaluc\'{\i}a grant  FQM-213, with FEDER funds.
}

\author[F.J.~Palomo]{Francisco J. Palomo${}^\diamond$}
%\curraddr{}
\email{fjpalomo@ctima.uma.es}
\thanks{${}^\diamond$ Supported by   the Spanish MEC grant MTM2016-78807-C2-2-P
and by the Junta de Andaluc\'{\i}a  grant  FQM-4496, with FEDER funds.}



\maketitle 

\begin{abstract}
The invariant metric affine connections on Berger spheres which are Einstein with skew torsion are determined in  both Riemannian and Lorentzian signature. Expressions of such connections are explicitly given. In particular, every Berger sphere with Lorentzian signature admits invariant metric affine connections which are Einstein with skew-torsion up to $\mathbb{S}^3$. For Riemannian signature, the existence of such connections strongly depends on the dimension of the sphere and on the scale of the deformation used for the Berger metric. In particular, there are Riemaniann Berger spheres, not Einstein, which admit invariant   Einstein with skew-torsion affine connections.
\end{abstract}

%In fact, for each odd dimension, there are many Riemaniann Berger spheres, not Einstein, which admit invariant metric affine connections Einstein with skew-torsion. MEJOR??

\section{Introduction}

Marcel Berger classified the simply connected normal homogeneous Riemannian manifolds with strictly positive sectional curvature  \cite{berger} (the Berger list was incomplete, missing one example, nowadays called the Aloff-Walach space $\mathfrak{W}_{1,1}$). In that paper, Berger found a family of $3$-dimensional Riemannian manifolds diffeomorphic to the $3$-sphere with strictly positive sectional curvature but not constant. These metrics can be obtained as a deformation of the round metric on $\mathbb{S}^{3}$ by changing the relative scale of the fibers of the Hopf fibration (see details in Section 2). These metrics, called now the \emph{Berger metrics}, provided counter-examples to several geometric conjectures. For example, Klingenberg proved %\cite[Chapter~3 and 5]{cheeger-ebin}  
that for an orientable
even-dimensional manifold $M$ whose sectional curvatures $K_{M}$ lie in the interval $[k_1, k_2]$ with $k_1>0$, every closed geodesic in $M$ has length at least $2\pi/\sqrt{k_2}$, while   Berger metrics show that such a result does not hold in general (see \cite[Chap. 3 and 5]{cheeger-ebin} for details and more references). Similar metrics for all  odd-dimensional spheres were constructed by I.~Chavel \cite{chavel}, and A.~Weinstein realized these metrics as distance spheres in the complex projective space \cite{weinstein}. On the other hand, the first nonzero eigenvalue of the Laplace operator of Berger spheres shows surprising properties which were first pointed out by H.~Urakawa \cite{urakawa}. Later, a nice description and extension of these properties from the submersion theory perspective was obtained in \cite{bergery}, where the reader will find more details and explanations about these facts.

From another point of view, the study of submanifolds immersed in Berger spheres has also been an active research field. Recall, for example, the works on Willmore surfaces \cite{barros} and minimal surfaces \cite{torralbo}. 
Moreover, for suitable scales of the fibers of the Hopf fibration, we can also consider the Berger spheres with Lorentzian signature (see Section 2 for details).   Lorentzian Berger spheres were studied in \cite{palomo}, where several properties of the conjugate points along their lightlike geodesics were obtained. From the point of view of dynamical systems, we would like to mention the recent interest in the study of magnetic curves in Berger spheres \cite{inoguchi}.

The Berger spheres (both Riemannian or Lorentzian signature) are  odd-dimensional spheres $\mathbb{S}^{2n+1}$ endowed with a $\SU(n+1)$-invariant metric $g_\varepsilon$ and so with constant sectional curvature but not Einstein, in general. Among the generalizations of the Einstein condition, this paper will focus on the notion of Einstein manifold with skew-torsion as was introduced in \cite{AgriFerr}. An Einstein manifold with skew-torsion is a triple $(M,g, \nabla)$ where $(M,g)$ is a semi-Riemannian manifold endowed with a metric  affine connection $\nabla$ which shares geodesics with the Levi-Civita connection and such that the symmetric part of the Ricci tensor of $\nabla$ is a multiple of the metric tensor $g$. In order  not to be so long, we recommend the complete and nice survey \cite{surveyagricola}, which explains the role of the  affine connections with torsion in Mathematics and Physics 
and provides a wide family of examples of connections with torsion in different contexts.
In this paper, we explicitly determine and describe the $\SU(n+1)$-invariant  affine connections $\nabla$ on the Berger spheres in such a way that $(\mathbb{S}^{2n+1}, g_{\varepsilon}, \nabla)$ is an Einstein manifold with skew-torsion. 
According to our conventions, for every $\varepsilon \neq 0$ the corresponding Berger sphere $(\mathbb{S}^{2n+1}, g_{\varepsilon})$ is a Riemannian manifold for $\varepsilon <0$ and Lorentzian for $\varepsilon >0$. 
This paper can be seen as a continuation of \cite{DraperGarvinPalomo}, which studies the case of the usual round  Riemannian metric of constant sectional curvature $1$, which corresponds to the value 
 $\varepsilon =-1$. It is proved in \cite{DraperGarvinPalomo} that, in this case, there are always Einstein connections with skew-torsion, although only for $n=3$ and $n=1$ there are nontrivial solutions, that is,  different from the Levi-Civita connection. 
Now the obtained results are qualitatively different: take into account that the Levi-Civita connection is not longer an Einstein connection with skew-torsion (for $\varepsilon\ne-1$), so any solution is nontrivial, and moreover, there are solutions for all the odd dimensions except for $\esf^3$.
In order to be precise, let us enunciate the main result.

\begin{theorem}\label{th_main}
Let $\esfera$ be the  odd-dimensional sphere viewed as homogeneous space $\SU(n+1)/\SU(n)$
and endowed with the invariant Berger metric $g_\varepsilon$ determined by 
Eq.~\eqref{eq_lag}. Take $\Phi,\eta,\psi,\Theta ,\widetilde{\Theta},\hat\psi$ the invariant tensors given by Eqs.~\eqref{sasaki},
\eqref{Lastetas} and \eqref{defhatpsi}. Let us denote by $\nabla^{g_{\varepsilon}}$ the Levi-Civita connection of $g_{\varepsilon}$.
Then, for the different values of $n$, we have
\begin{enumerate}
%\item[\boxed{$n\ne 2,3$}]
\item[\boxed{ {n\ne  2,3}}] 
The invariant metric connections with skew-torsion are:
$$
\nabla_{X}Y=\nabla^{g_{\varepsilon}}_{X}Y+s\,\big(\Phi(X,Y)\,\xi+\varepsilon(\eta(X)\psi(Y)-\eta(Y)\psi(X))\big)
$$
for any $s\in\RR$.
\begin{itemize}\item
If $n=1$,  they are      Einstein with skew-torsion affine connections  if and only if $\varepsilon=-1$.
In this case, every choice of the parameter $s$ gives such a connection.
\item
If $n>3$, they are   Einstein with skew-torsion affine connections   if and only if
$$
s^2= {\left(\frac{n+1}{n-1}\right)\,   \left(\frac{\varepsilon+1}\varepsilon\right)  }.
$$
\end{itemize}
\item[\boxed{ {n=3}}] 
The invariant metric connections with skew-torsion are:
$$
\begin{array}{ll}\vspace{1pt}
\nabla_{X}Y =\nabla^{g_{\varepsilon}}_{X}Y&+s\,\big(\Phi(X,Y)\,\xi+\varepsilon(\eta(X)\psi(Y)-\eta(Y)\psi(X))\big)\\
&+s_1 \, {\Theta}(X,Y)+s_2\,\widetilde{\Theta}(X,Y)
\end{array}
$$
for any $s,s_1,s_2\in\RR$.
They are     Einstein with skew-torsion affine connections if and only if
$$
\varepsilon{s^2} + {s_1^2+s_2^2} =2(\varepsilon+1).
$$
\item[\boxed{ {n=2}}] 
The invariant metric connections with skew-torsion are:
$$
\begin{array}{ll}\vspace{1pt}
\nabla_{X}Y =\nabla^{g_{\varepsilon}}_{X}Y&+s\,\big(\Phi(X,Y)\,\xi+\varepsilon(\eta(X)\psi(Y)-\eta(Y)\psi(X))\big)\\
\vspace{1pt}
&+s_3\,\big(\Phi(\hat\psi(X),Y)\,\xi+\varepsilon(\eta(X)\psi(\hat\psi(Y))-\eta(Y)\psi(\hat\psi(X)))\big)\\
\vspace{1pt}
&+s_4\,\big(- g_\varepsilon( \hat\psi(X),Y)\,\xi+\varepsilon(\eta(X)\hat\psi(Y)-\eta(Y)\hat\psi(X)) \big)
\end{array}
$$ 
for any $s,s_3,s_4\in\RR$.
They are      Einstein with skew-torsion affine connections  if and only if
$$
 {s^2+s_3^2+s_4^2}=3\,\frac{\varepsilon+1}{\varepsilon}.
$$
\end{enumerate}
\end{theorem}

In particular, we recover the above mentioned well-known result  (see \cite[9.81]{besse} for a completely different proof), asserting that 
the unique Einstein Riemannian metric in the family $g_{\varepsilon}$ with $\varepsilon < 0$ is the canonical round metric $g_{-1}$ of constant sectional curvature $1$. 
We also note that, up to $\mathbb{S}^{5}$, the contraction of the $3$-differential forms corresponding with the torsions 
(see   Section~2) are exactly the closed $2$-forms used in \cite{inoguchi} to define magnetic curves.

\vspace{3pt}
%In particular, Levi-Civita connection is Einstein  if and only if $\varepsilon=-1$. 
As a direct corollary of Theorem~\ref{th_main},  we can analyze the existence of    Einstein with skew-torsion affine connections according to the values of $\varepsilon$ and $n$ (see Table~1).\vspace{1pt}

\begin{quote}

For Lorentzian signature ($\varepsilon >0$), the Berger sphere $(\mathbb{S}^{2n+1}, g_{\varepsilon})$ always admits a 
$\SU(n+1)$-invariant   Einstein with skew-torsion affine connection up to $n=1$. For $n\in \{2,3\}$, the set of such connections is parametrized by the points of an ellipsoid while reduces to just two connections for $n\geq 4$.
\end{quote}\vspace{1pt}

%This is already interesting, since the Levi-Civita connection of any  invariant Lorentzian metric is never Einstein with skew-torsion.

This is already interesting, since Lorentzian metrics $g_\varepsilon$ are never Einstein (in the usual sense).
Thus, there would be another natural choice for a distinguished affine connection in the Lorentzian odd-dimensional spheres of dimension at least $9$.\vspace{1pt}

\begin{quote}
For Riemannian signature ($\varepsilon < 0$), the Berger sphere $(\mathbb{S}^{3}, g_{\varepsilon})$ admits   
$\SU(2)$-invariant Einstein with skew-torsion  affine connections  only for $\varepsilon =-1$, and the set of such connections is parametrized by a line. 
For $(\mathbb{S}^{5}, g_{\varepsilon})$, the situation depends on $\varepsilon$; for $\varepsilon > -1$ there is no such connection, for $\varepsilon =-1$ there is only one and for $\varepsilon < -1$ the set of such connections is parametrized by an ellipsoid. 
For $(\mathbb{S}^{7}, g_{\varepsilon})$, the set of 
$\SU(4)$-invariant  Einstein with skew-torsion  affine connections is parametrized by a one-sheet hyperboloid
 if $\varepsilon > -1$, by a cone if $\varepsilon =-1$   and by a    two-sheets hyperboloid if $\varepsilon < -1$. 
Finally, for $n>3$, $(\mathbb{S}^{2n+1}, g_{\varepsilon})$ admits a $\SU(n+1)$-invariant Einstein with skew-torsion  affine connection  only for $\varepsilon \leq -1$; there is only one if $\varepsilon =-1$ and exactly two if $\varepsilon < -1$.
\end{quote}\vspace{1pt}

The sphere $\mathbb{S}^7$ is exceptional (among other properties) from this point of view. In fact,  $\mathbb{S}^7$  is the unique sphere such that all Riemannian metrics $g_{\varepsilon}$ admit invariant   Einstein with skew-torsion affine connections. In the opposite position we find $\mathbb{S}^3$, which admits invariant   Einstein with skew-torsion affine connections only for the round metric. 
%It is also remarkable that the round sphere is not representative at all among the Berger Riemannian spheres.

In our view, we can divide the Riemannian Berger metrics with similar behaviors into three groups: the round metric (alone), those metrics with distort enlarging the vertical direction, and those ones shortening the vertical direction.
 
\vspace{0.2cm}

%Recall now that magnnetic field on a  semi-Riemannian manifold $(M,g)$ is a closed $2$-differential form $\omega$. A curve $\gamma: I \to M$ is said to be a magnetic curve when $\gamma$ satisifies the Lorentz equation
%$
%\nabla^{g}_{\gamma' }\gamma'=F(\gamma' ),
%$
%where $g(X,Y)=g(F(X),Y)$ for $X,Y \in \mathfrak{X}(M)$
%Every  odd-dimensional sphere admits a Sasakian structure (see Section 2) and several tensors closely related with such structure has been used in theorem \ref{th_main}. 
%The study of magnetic curves in Berger spheres (Inoguchi)

 \medskip

 The structure of the paper is the following:  
 We recall in Section~2 the algebraic tools to work with invariant affine connections in a reductive homogeneous space, mainly Nomizu's Theorem.
 We also describe the sphere as $\SU(n+1)$-homogeneous space and list the invariant tensors which will need for the concrete descriptions of the connections in Theorem~\ref{th_main}. As this paper continues the work in \cite{DraperGarvinPalomo}, we have repeated only the necessary material to achieve and write down the results.
 Section~3   contains the proof of Theorem~\ref{th_main} in a schematic   way.

 %: most tedious computations are omitted but the results can be checked easily by the reader.  

%Organization and results of this paper....
% 
%\begin{coro}\label{co_main}
%Consider  the Berger spheres $(M,g_\varepsilon)$ as in the previous theorem.
%
%\begin{itemize}
%\item If $n=1$, there are Einstein connections only for the round sphere ($\varepsilon=-1$), and in this case 
%the set of Einstein connections is parametrized by a line.
%\item If $n=2$, there are no Einstein connections if $\varepsilon>-1$, there is only one if $\varepsilon=-1$, and the set 
%of Einstein connections is parametrized by an ellipsoid
%if $\varepsilon<-1$.
%\item If $n=3$, the set of of Einstein connections is parametrized by a one sheet hyperbolloid if $\varepsilon>-1$, by a cone if $\varepsilon=-1$ and by a two sheets hyperbolloid   if $\varepsilon<-1$.
%\item  If $n>3$, there are no Einstein connections if $\varepsilon>-1$, there is only one if $\varepsilon=-1$ and there are just two if $\varepsilon<-1$.
%\end{itemize} 
%\end{coro}

%%%%%%%%%%%%%%%%%%%%%%%%%%%%%%%%%%%%%%%%%%%%%%%%%%%%%%%%%%%%%%%%%%%%%%%%%%%%%%%%%%%%%%%%%%%%%%%%%%%%%%%%%%

 \section{Preliminaries.} 
 
 %We recall from ... several key facts on the study of these spheres.
 \subsection{The odd-dimensional spheres as homogeneous spaces for the special unitary group and algebraic tools.}
 
 Throughout this paper we consider the sphere $M=\esfera\subset \mathbb C^{n+1}$ as homogeneous space under the action of the special unitary group $\SU(n+1)$. Invariant affine connections on $M$ have been studied on \cite{DraperGarvinPalomo}, where those ones compatible with the round (Riemannian) metric have been explicitly described, and more precisely,  being besides Einstein with skew-torsion in the sense of \cite{AgriFerr}. Part of the results of \cite{DraperGarvinPalomo} are the starting point now. We briefly recall the situation.\vspace{2pt}
 
 The special unitary group 
 $
 G=\SU(n+1)=\{ A\in \mathcal M _{n+1}(\mathbb C) :  A\bar A^t=I_{n+1} , \, \det(A)=1\}
$
 acts smoothly on 
 $$
 M=\esf^{2n+1}=\{ (z_1,\cdots , z_{n+1})^t\in \mathbb C^{n+1} :  \sum _{i=1}^{n+1}z_i\bar z_i=1\} ,
 $$
 by matrix multiplication. As the action is transitive, and the isotropy group of the point 
 $o=(0,\cdots ,0,1)^t \in M$ is $H_o=\{\tiny\left(
\begin{array}{c|c}
  B&0    \\ \hline
  0& 1  
\end{array}
\right):B\in\SU(n)\}$, %then 
 the sphere $\esfera$ is diffeomorphic to the homogeneous space $G/H=\frac{\SU(n+1)}{\SU(n)}$. Moreover, this homogeneous space is \emph{reductive}, that is, if $\g$ is the Lie algebra of $G$ and $\hh$ the Lie algebra of (the connected group)
 $H$, there is a complementary subspace $\mm$ of $\hh$ such that $[\hh,\mm]\subset\mm$. 
 There is always a natural identification of $\mm$ with the tangent space at  $o=\pi(e)$, provided by the differential of the projection $\pi\colon G\to M=G/H$, which gives a linear isomorphism
$(\pi_{*})_e\vert_{\mathfrak{m}}\colon\mathfrak{m}\to T_{o}M$.
 For our precise example, 
 $\mathfrak g=\mathfrak{su}(n+1)$, 
 $
\mathfrak h=\{\tiny
\left(
\begin{array}{c|c}
  B&0    \\ \hline
  0& 0  
\end{array}
\right): B\in \mathfrak{su}(n)
\}
$ and the complementary 
$\mm$ is unique, namely, 
$$
\mathfrak m=\left\{
\left(
\begin{array}{c|c}
  -\frac{a}{n}I_n&z    \\ \hline
  -\bar z^t& a  
\end{array}
\right) \in \mathcal M_{n+1}(\mathbb C) :  z\in \mathbb C^n, a\in \mathbb R\mathbf{i}
\right\}.
$$
Under the natural identifications $\hh\cong\mathfrak{su}(n)$ and $\mm\cong \mathbb C^n\oplus \mathbb R\mathbf{i}$, 
 the action $[\hh,\mm]\subset\mm$ becomes $B\cdot (z,a)=(Bz,0)$.\medskip

%\textcolor{violet}{ The study of invariant affine connections in reductive  homogeneous spaces   can be translated to an algebraic setting thanks to the Nomizu's Theorem (cf. \cite{teoNomizu}). In these terms, it is not difficult to compute connections with suitable properties, and try to choose a convenient connection in some sense.} \margen{este color es que cambiaré el contenido}

%%DEFINICION DE INVARIANTE LA QUITO 
%Let $G$ be a Lie group acting transitively on $M$. We write a dot to denote the action of $G$ on $M$ and so, for $\sigma \in G$, the left translation by $\sigma$ will be given by $\tau_{\sigma}(p)=\sigma \cdot p$ for all $p\in M$.
%For each $\sigma \in G$ and $X\in \mathfrak{X}(M)$, the vector field $\tau_{\sigma}(X)\in \mathfrak{X}(M)$ is given at each $p\in M$ by
%\[
% (\tau_{\sigma}(X) )_p:=(\tau_{\sigma})_{*}(X_{\sigma^{-1}\cdot p}).
%\]
%An affine connection $\nabla$ on $M$ is said to be $G$-\textit{invariant}  if, for each $\sigma \in G$ and for all $X,Y\in\mathfrak{X}(M)$,  \margen{copiado de version previa de 3-sasaski}
%$$
%\tau_{\sigma}(\nabla_{X}Y)=\nabla_{_{\tau_{\sigma}(X)}}\tau_{\sigma}(Y).
%$$ 

The study of invariant affine connections in reductive  homogeneous spaces   can be translated to an algebraic setting thanks to   Nomizu's Theorem (cf. \cite{teoNomizu}).

%Let $H=\{\sigma \in G:\sigma \cdot o=o\}$ be the isotropy subgroup of a fixed point $o\in M$, so that there exists a diffeomorphism between $M$ and $G/H$. If $H$ is connected, the homogeneous space $M=G/H$ is said to be \textit{reductive} when the Lie algebra $\mathfrak{g}$ of $G$ admits a vector space decomposition
%\begin{equation}\label{eq_descomposicionreductiva}
%\mathfrak{g}=\mathfrak{h}\oplus\mathfrak{m},
%\end{equation}
%for $\mathfrak{h}$ the Lie algebra of $H$ and $\mathfrak{m}$ a complementary subspace  such that $[\mathfrak{h},\mathfrak{m}]\subset\mathfrak{m}$. In this case, $\mathfrak{g}=\mathfrak{h}\oplus\mathfrak{m}$ is called a \emph{reductive decomposition} of $\mathfrak{g}$.  The differential map $\pi_{*}$ of the projection $\pi\colon G\to M=G/H$ gives a linear isomorphism
%$(\pi_{*})_e\vert_{\mathfrak{m}}\colon\mathfrak{m}\to T_{o}M$, where $o=\pi(e)$. Nomizu's Theorem   \cite{teoNomizu} can be formulated as follows: \margen{caso conexo}

%\renewcommand{\thetheorem}{\Alph{theorem}}
%\setcounter{theorem}{13}

\begin{theorem}\label{nomizu}
Let $G/H$ be a reductive homogeneous space %($H$ connected) 
with a fixed reductive decomposition $\mathfrak{g}=\mathfrak{h}\oplus\mathfrak{m}$ and $H$ connected.
Then, there is a bijective correspondence between the set of   $G$-invariant affine connections  on $G/H$ and the vector space of bilinear maps $\alpha\colon \mathfrak{m}\times\mathfrak{m}\to\mathfrak{m}$ such that 
\begin{equation*}\label{casoconexo}
[h,\alpha(X,Y)]=\alpha([h,X],Y)+\alpha(X,[h,Y])
\end{equation*}
for all $X,Y\in\mathfrak{m}$ and $h\in\mathfrak{h}$. 

Moreover, given a $G$-invariant affine connection $\nabla$ 
%compatible with the metric $g$ and having  skew-symmetric torsion $T^{\nabla}$, 
and $\alpha $ the associated bilinear map, the torsion and curvature tensors of 
%the $G$-invariant affine connection 
$\nabla$
are given by   
\begin{align}\label{tor}
T_\alpha(X,Y)&= \alpha  (X,Y)-\alpha  (Y,X)-[X,Y]_\mathfrak{m}, \\
\label{cur}
 R_\alpha(X,Y,Z)&=
 \alpha  (X,\alpha (Y, Z))-\alpha  (Y,\alpha  (X,Z))-\alpha  ( [X,Y]_{\mathfrak{m}}, Z)-[[X,Y]_{\mathfrak{h}} , Z],
\end{align}
for any $X,Y,Z\in \mathfrak{m}$, where   $[\ ,\ ]_{\mathfrak{h}}$ and   $[\ ,\ ]_{\mathfrak{m}}$ denote the composition of the bracket
$ [\mathfrak{m},\mathfrak{m}]\subset \mathfrak{g}$ with    the projections 
%$\pi_{\mathfrak{h}}$ and $\pi_{\mathfrak{m}}$  
of $\mathfrak{g}= \mathfrak {h}\oplus \mathfrak{m}$ 
on each factor. 
\end{theorem}

%As $\mm\cong T_oM$, 
These expressions give the torsion $T $  and the curvature $R $ at the point $o$ (by means of the above identification of $\mathfrak{m}$ with $T_oM$ via $\pi_{*}$). The invariance permits to recover the whole tensors.  \smallskip

In particular, there is also a bijective correspondence among the set of   $G$-invariant affine connections  on $G/H$
and the vector space of   homomorphisms of $\mathfrak{h}$-modules 
$\mathrm{Hom}_{\mathfrak{h}}(\mathfrak{m}\otimes\mathfrak{m},\mathfrak{m})$, where the morphism related to $\alpha$ is given by $\mathfrak{m}\otimes\mathfrak{m}\to\mathfrak{m}$, $X\otimes Y\mapsto \alpha (X,Y)$.

The   vector space  $\mathrm{Hom}_{\mathfrak{h}}(\mathfrak{m}\otimes\mathfrak{m},\mathfrak{m})$ has been studied in 
\cite{DraperGarvinPalomo}
for our spheres. Concrete bases are exhibited in  \cite[Eqs. (25), (31), \S6.1 and \S7]{DraperGarvinPalomo}, and the dimension turns out to be $7$, $9$, $13$ and $27$   for $n\ge4$, $n=3$, $n=2$ and $n=1$ respectively. 
\medskip

If a reductive homogeneous space $G/H$ is endowed with a $G$-invariant metric $g$,
%\footnote{Any such metric is determined by the nondegenerate symmetric bilinear map $g_o\colon T_oM\times T_oM\to\mathbb R$ by means of the identification of $\mathfrak{m}$ with $T_oM$ via $\pi_{*}$. We also call $g$ for referring $g_o$.}
 there seems natural to ask for the invariant connections to be also compatible with the metric  ($\nabla g=0$). 
In such a case, the affine connection  is determined by the torsion.  Nomizu's theorem can be also adapted to this setting. For instance from \cite[Theorem 2.7, Remark 2.8]{DraperGarvinPalomo}, if we also denote  by
$g\colon \mathfrak{m}\times\mathfrak{m}\to\RR$   the $\hh$-invariant nondegenerate symmetric bilinear map induced by $g\vert_{T_oM}$,   %by means of the identification of $\mathfrak{m}$ with $T_oM$ via $\pi_{*}$.
then we have,
 
\begin{lemma} 
A $G$-invariant affine connection   is metric  if and only if 
$
 % g  (\alpha_{_{\nabla}}(C,A),  B)  + g (A,\alpha_{_{\nabla}} (C,B))=0
  \alpha (X,-)\in\mathfrak{so}(\mathfrak{m},g)$  for all $X\in\mathfrak{m}$. 

  %for any $A,B,C \in \mathfrak{m}$, 
%  for all $X\in\mathfrak{m}$, 
%where 
 \end{lemma}
 
  In other words, there is a bijective correspondence between the set of $G$-invariant affine connections   compatible with the invariant metric $g$ on $M=G/H$ and the vector space $
\hom_{\mathfrak{h}}(\mathfrak{m}, \wedge^2\mm)$, taking into account that $\wedge^2\mm$ and $\sof(\mm,g)$ are isomorphic $\hh$-modules.
% The affine connection in this case is determined by the torsion. %, being $\nabla=\nabla^g+\frac12T$ SOLO para antisimetrico
% 
 % 
  In particular, the set of $G$-invariant metric affine connections    is parametrized by a number of parameters which does not depend on the concrete choice of the invariant metric. \smallskip
  
  Let us return to the odd-dimensional spheres $\esfera=\frac{\SU(n+1)}{\SU(n)}$.
  The Berger metric $g_\varepsilon$ for each $0\ne\varepsilon\in\mathbb R$ can be introduced as follows.
  The bilinear map
  \begin{equation}\label{eq_lag}
  g_\varepsilon\colon\mm\times\mm\to\mathbb R,\qquad  g_\varepsilon((z,a),(w,b)):=\Re(z^t\bar w)+\varepsilon ab
 \end{equation}
is an $\hh$-invariant nondegenerate symmetric bilinear map, so that  it determines 
%(after identifying with  $(g_\varepsilon)_o\colon T_oM\times T_oM\to\mathbb R$) 
an invariant metric $g_\varepsilon\in\mathcal T^{0,2}(M)$.  The round Riemannian metric used in 
\cite{DraperGarvinPalomo} corresponds with $\varepsilon=-1$. In general, when $\varepsilon<0$, the obtained metric 
$g_\varepsilon$ is Riemannian, while $g_\varepsilon$ is Lorentzian for $\varepsilon>0$. Indeed, $g_\varepsilon\vert_{\mathbb C^{n}}$ is positive definite, independently of $\varepsilon$, and the index is determined by $g_\varepsilon((0,\mathbf{i}),(0,\mathbf{i}))=-\varepsilon$. Moreover this family $\{g_\varepsilon: \varepsilon\ne0\}$ exhausts the set of all invariant metrics up to nonzero scalar multiple, as explained in \cite[Remark~3.1]{DraperGarvinPalomo}.

Now we can use, for any of the metrics $g_\varepsilon$, the results in \cite[Prop.~4.4, Prop.~5.3, Prop.~6.3, \S7]{DraperGarvinPalomo}, %which assert that 
according to which
$$
\dim\hom_{\mathfrak{su}(n)}(\mathfrak{m}, \wedge^2\mm)=3,5,7,9
$$
respectively for $n\ge 4$, $n=3$, $n=2$ and $n=1$.
%About our precise computations in the next section, 
Anyway,
we will find a concrete description of the metric connections without using the previous knowledge of the dimension of the vector space.
%this result, but coherent with it.
\medskip

%\textcolor{blue}{definir $\omega_\alpha$. En tercer lugar, Now, we can also translate the problem of la torsión antisimétrica a un contexto algebraico.}

We want to find connections as good as possible. At this point we recall that a  semi-Riemannian manifold $(M,g)$ endowed with a metric affine connection $\nabla$  is called
\emph{Einstein with skew-torsion} (see \cite{AgriFerr} for a variational justification of this notion) if:
%\begin{itemize}
%\item[i)] 
$\nabla$  has totally skew-symmetric torsion, that is, 
$
\omega(X,Y,Z):=g(T(X,Y),Z)
$
is a 3-form on $M$ for $T$ the usual torsion tensor;
%\item[ii)]  
and the symmetric part of the Ricci tensor satisfies
\begin{equation}\label{eq_Einsteincondition}
 %$
\mathrm{Sym}(\mathrm{Ric}^{\nabla})=\frac{s^{\nabla}}{\dim M}\, g,
%$
\end{equation} 
where   $s^{\nabla}$ denotes the scalar curvature.
We will also say that $\nabla$ is an Einstein with skew-torsion connection.
%\end{itemize} %\margen{comentar que para s cte equivale a ser multiplo de g}
%In our situation of constant scalar curvature, condition ii) is equivalent to the symmetrized tensor of the Ricci tensor to be a scalar multiple of the metric. \margen{y si no ¿también?}
In the above setting, the condition   of having skew-symmetric torsion is   a linear condition on the parameters of the vector space
$\hom_{\mathfrak{su}(n)}(\mathfrak{m}, \wedge^2\mm)$, and the obtained maps can be proved to be in one-to-one correspondence with $\mathrm{Hom}_{\mathfrak{h}}(  \wedge^3 \mathfrak{m},\RR)$. On the contrary,   
Eq.~\eqref{eq_Einsteincondition} 
%the \emph{Einstein condition} ii)  
%(equivalent to the symmetrized tensor of the Ricci tensor to be a scalar multiple of the metric) 
provides a set of quadratic equations on these parameters, so that the set of solutions could be   empty.\smallskip

%\begin{lemma}  \label{igualdadconjuntos} Let $(M=G/H,g)$ be a reductive homogeneous space with a $G$-invariant metric $g$. The set ............. is bijective to $\mathrm{Hom}_{\mathfrak{h}}(\mathfrak{m} \wedge \mathfrak{m} \wedge \mathfrak{m},\RR)$. .................
%\end{lemma}
%
%\textcolor{blue}{Finalmente, estaremos interesados  en hallar las conexiones cuyo tensor de Riccis sea ..... Obviamente esa condición ya no es lineal, pero dará lugar a un conjunto de ecuaciones cuadráticas.....}

In what follows,  we occasionally write connections, metric connections and metric connections with skew-torsion for referring to the bilinear maps related to those connections through Nomizu's Theorem.

%\textcolor{blue}{por abuso de notacion hablare de conexiones para sus apl bil asociadas. decir en preliminares}

%%%%%%%%%%%%%%%%%%%%%%%%%%%%%%%%%%%%%%%%%%%%%%%%%%%%%
\subsection{Invariant tensors}

We will derive concrete expressions for the $\SU(n+1)$-invariant affine  connections using some invariant tensors on $M=\esf^{2n+1}$. For most of the cases ($n\ne2,3$),   the   connections will be  described using only  the usual Sasakian structure on $\esf^{2n+1}$. We recall briefly here the involved tensors (see \cite{Blair}).
For $X\in \mathfrak{X}(\esf^{2n+1})$, the decomposition of $\mathbf{i}  X$ in tangent and normal components defines the $(1,1)$-tensor field 
  $\psi$ and the 1-differential form $\eta$ as follows,
  \begin{equation}\label{sasaki}
\mathbf{i}  X=\psi(X)+\eta(X)N,
\end{equation}
where $N$ denotes the unit outward normal vector field to $\esf^{2n+1}$.
Let $\xi\in \mathfrak{X}(\esf^{2n+1})$ be the Killing vector field (for all $g_\varepsilon$) defined by
$
\xi_{z}=-\mathbf{i}z
$
at any $z\in \esf^{2n+1}$.
%
%
%
%Thus, if we denote by $\nabla^g$ the Levi-Civita connection of $g$, the following properties hold
%\begin{equation}\label{eq_propiedadesSasakiana}
%\begin{array}{ll}
%\eta(X)=g(X,\xi),\qquad  &g(\psi(X),\psi(Y))=g(X,Y)-\eta(X)\eta(Y),
%\\
%\psi^2=-\mathrm{Id}+\eta \otimes \xi,\quad\qquad&(\nabla^{g}_{X}\psi)Y=g(X,Y)\xi-\eta(Y)X,
%\end{array}
%\end{equation}
%for any $X,Y\in \mathfrak{X}(\SS^{2n+1})$. These conditions imply several further relations, for instance we also have that $\eta \circ \psi=0$ and $\psi(\xi)=0$. 
The Sasakian form   $\Phi\in\Omega^2(\esfera)$ is defined by
$\Phi(X,Y)=g_\varepsilon(X, \psi(Y))$ (note that $\Phi$ is independent of $\varepsilon$). 
%Moreover, the vector field $\xi$ is Killing (i.e., $g(\nabla^{g}_{X}\xi,Y )+g(\nabla^{g}_{Y}\xi,X )=0$)
%% for all $X,Y\in\SS^{2n+1}$)  
%and therefore $\nabla^{g}_{X}\xi=-\psi(X)$ and $2\Phi =d\eta$.

The tensor fields $\psi,\eta$ and $\Phi$ are $\SU(n+1)$-invariant \cite[Lemma~3.2]{DraperGarvinPalomo}, so that they are determined by the value at the point $o\in M$. Note that, if 
$X=(z,a),Y=(w,b), Z=(u,c)\in \mathfrak{m}=T_oM\cong \mathbb C^n\oplus \mathbb R\mathbf{i}$, then
\begin{equation}\label{eq_tensoresenm}
\begin{array}{ll}
\psi(X)=(\mathbf{i} z,0)\in\mm,\qquad&
\eta(X)=\mathbf{i} a\in\RR, \\
\xi=(0,-\mathbf{i})\in\mm,&
\Phi(X,Y)=g_{\varepsilon}(X,\psi(Y))=-\mathrm{Im}( \,\overline{z}^tw)\in\RR.
\end{array}
\end{equation}\medskip

 In the case $n=3$, the geometric structure of the sphere $\esf^7$ is richer. In fact, $\esf^7$ is a $3$-Sasakian manifold as follows. If we realize $\esf^{7}$  as a hypersurface of $\mathbb H^2$ by means of
\begin{equation*} 
\esf^{7}\subset \CC^{4}\cong\mathbb H^{2}, \quad z=(z_{_1},z_{_2},z_{_3},z_{_4})\mapsto (z_{_1}+z_{_2}\,{\bf j},z_{_3}+z_{_4}\,{\bf j}),
\end{equation*}
we can consider the vector fields
$\xi_{1},\xi_{2},\xi_{3}\in \mathfrak{X}(\esf^{7})$ given by $\xi_{1}(z)=-{\bf i}z$, $\xi_{2}(z)=-{\bf j}z$ and
$\xi_{3}(z)=-{\bf k}z$, respectively. For any vector field $X\in \mathfrak{X}(\esf^{7})$, the decompositions of $\mathbf{i}\,X$, $\mathbf{j}\,X$ and $\mathbf{k}\,X$ into tangent and normal components  determine as above the $(1,1)$-tensor fields $\psi_{1},\psi_{2}$ and $\psi_3$ and the differential $1$-forms $\eta_{1},\eta_{2}$ and $\eta_{3}$ on $\esf^{7}$, respectively.
 The
 three  Sasakian structures $(\xi_{i},\eta_{i}, \psi_{i})$ are compatible in the sense that (indices modulo 3)
 $$
 %\begin{array}{l}
[\xi_{i},\xi_{i+1}]=2\, \xi_{i+2},\quad
\psi_{i+1}\circ \psi_{i}=-\psi_{i+2}+\eta_{i}\otimes \xi_{i+1},\quad
\psi_{i}\circ \psi_{i+1}=\psi_{i+2}+\eta_{i+1}\otimes \xi_{i}.
%\end{array}
$$
 The Sasakian structure $(\xi_{1},\eta_{1}, \psi_{1})$ is $\SU(4)$-invariant, but this is not the case for the other two indices. Fortunately, it is found in \cite[Proposition~5.9]{DraperGarvinPalomo} the following $\SU(4)$-invariant tensor
  on $\esf^{7}$, the
 $3$-differential form  
 $
 %\begin{equation}\label{LaOmega}
 \Omega=%\frac{1}{2}(\eta_{2}\wedge d\eta_{2}-\eta_{3}\wedge d\eta_{3})=
  \eta_{2}\wedge \Phi_{2}-\eta_{3}\wedge \Phi_{3}.
 %\end{equation}
 $
 Hence the 2-forms  $\Theta ,\widetilde{\Theta}\in\Omega^2(\esf^7) $   defined by contraction as
\begin{equation}\label{Lastetas}
g_{-1}(\Theta (X,Y),Z)= \Omega(X,Y,Z)
,\qquad
g_{-1}(\widetilde{\Theta} (X,Y),Z)= \Omega(X,Y,\psi_{1}(Z))
\end{equation}
are also invariant.
Under the identification $\pi_{*}\colon\mathfrak{m}\to T_{o}\esf^7$,
\begin{equation*}%\label{eq_tensoresenm3}
\Omega(X,Y,Z)=-\Re(\det(z,w,u)),\qquad \Omega(X,Y,\psi_{1}(Z))=\mathrm{Im}(\det(z,w,u)); 
\end{equation*}  %%%cpbado, es un +
  that gives
\begin{equation}\label{eq_tensoresenm3}
 {\Theta} (X,Y)=-(\bar z\times\bar w,0),\qquad \widetilde{\Theta} (X,Y)=(\mathbf{i}\,(\bar z\times\bar w),0).
\end{equation} %%También comprobado
\smallskip

 Finally the sphere $\esf^5$ has   a distinguished invariant tensor too.   
Let us consider the sphere $\esf^6$ in $\RR^7$ with unit outward normal vector field $N$, and the almost complex structure $J$ on $\esf^6$   defined by $J(X)= N\times X$ for every $X\in \mathfrak{X}(\esf^6)$, for $\times$ a cross product in $\RR^7$. 
Thus, $(\esf^6, g_{-1} , J)$ is a nearly K\"{a}hler manifold \cite{Blair}. 
%Think of $\CC$ as $\RR\oplus\RR e_7$, and consider the
%totally geodesic embedding $\esf^{5}\to \esf^6$ defined by
%$$
%(z_{1},z_{2},z_{3})\in \esf^5\subset \CC^3\mapsto z_{1}  e_{1}+z_{2}  e_{2}+z_{3}  e_{4}+0e_7\in \esf^6 \subset \mathrm{Im}(\mathbb{O})\cong\mathbb{R}^7,
%$$
%%(with $e_7e_1=e_3  $, $e_7e_2=e_6  $, $e_7e_4=e_5  $), 
%which has unit normal vector field $\nu=-\frac{\partial}{\partial x_{7}}$.
%Here the basis is written so that $e_1\times e_2=\frac12 e_4$ and $e_7e_i=\frac12e_7\times e_i$ if $i\ne7$.
The sphere $\esf^5$ can be totally geodesic embedded in $\esf^6$ at the equator ($x_7=0$).
For every $X\in \mathfrak{X}(\esf^5)$, the decomposition of $J(X)$ in tangent and normal components with respect to the above isometric embedding $\esf^5 \subset \esf^6$ determines the $(1,1)$-tensor field $\widehat{\psi}$ and (again) the $1$-differential form $\eta$ on $\esf^5$ such that     
\begin{equation}\label{defhatpsi}
J(X)=\widehat{\psi}(X)+\eta(X)\nu,
\end{equation}
where $\nu$ is the unit normal vector field along the above totally geodesic embedding \cite{Blair}.
%Thus, $\xi$, $\eta$ and $\widehat{\psi}$ form an 
 The almost contact metric structure $(\xi,\eta,\widehat{\psi})$ on $\esf^5$ is $\SU(3)$-invariant (but not Sasakian).  
The linear map  $\widehat{\psi}_{o}\colon T_{o}\esf^5\to T_{o}\esf^5$ corresponds under the identification $\pi_{*}\colon\mathfrak{m}\to T_{o}\esf^5$ with the endomorphism of $\mathfrak{m}\cong\CC^2\oplus\RR\mathbf{i}$ given by 
\begin{equation}\label{eq_tensoresenm2}
\widehat{\psi}(z,a)=(\theta(z),0),
\end{equation}
for
$\theta\colon\CC^2\to\CC^2$ the $\hh=\mathfrak{su}(2)$-homomorphism given by
  $
  \theta\left(\begin{array}{l}  z_1\\   z_2\end{array}\right)=\left(\begin{array}{c}-\bar z_2\\ \bar z_1\end{array}\right)$
  \cite{DraperGarvinPalomo}.

%%%%%%%%%%%%%%%%%%%%%%%%%%%%%%%%%%%%%%%%%%%%%%%%%%%%%
%%%%%%%%%%%%%%%%%%%%%%%%%%%%%%%%%%%%%%%%%%%%%%%%%%%%%
 \section{Proof of the Main Theorem}
 
 For brevity, we take $X=(z,a)$, $Y=(w,b)$ and $Z=(u,c)$ arbitrary elements in $ \mathfrak{m}$.

\subsection{Case $\esf^{2n+1}\cong\frac{\SU(n+1)}{\SU(n)}$ for $n\ge4$.}

Recall from \cite[Proposition~4.3]{DraperGarvinPalomo} that an $\RR$-bilinear map $\alpha\colon\mm\times\mm\to\mm$ is $\hh$-invariant if and only if there exist $q_{1},q_{2},q_{3}\in \mathbb{C}$ and $t\in \mathbb{R}$ such that, for any $X,Y\in\mm$,
\begin{equation*}\label{eq_todoscason}
\alpha(X,Y)=\left( q_{1}bz+q_{2}aw,\, \mathbf{i}\, \left(tab+\mathrm{Im}(q_{3}\overline{z}^tw)\right)\right).
\end{equation*}

We easily check that  $\alpha\in\mathfrak{so}(\mathfrak{m},g_\varepsilon)$ when $t=\mathrm{Im}(q_{2})=0$ and $q_{1}=\varepsilon\overline{q_{3}}$, so that, by renaming, there are $q\in \mathbb{C}$ and $t\in \mathbb{R}$  such that 
\begin{equation*}\label{eq_metriccason}
\alpha(X,Y)=\left( -\varepsilon q\,bz+t\,aw,\, -\mathbf{i}\,  \mathrm{Im}(\bar q\,\overline{z}^tw)\right).
\end{equation*}
 The torsion tensor as in Eq.~\eqref{tor}  is
\begin{equation}\label{eq_torsioncason}
T{_{\alpha}}(X,Y)=\left(\Big(-\varepsilon q-t-\frac{n+1}{n}\Big)(bz-aw), (\mathrm{Re}(q)-1)(\overline{w}^tz-\overline{z}^tw)\right),
\end{equation}
which vanishes for $q=1, t=-\varepsilon -\frac{n+1}{n}$. This means that the Levi-Civita connection of $g_\varepsilon$ corresponds to
\begin{equation}\label{eq_LCcason}
\alpha_{g_\varepsilon  }(X,Y)=\left( -\varepsilon  \,bz-(\varepsilon +\frac{n+1}{n})\,aw,\,
 -\mathbf{i}\,  \mathrm{Im}( \,\overline{z}^tw)\right).
\end{equation}
If we compute the $(0,3)$-tensor $\omega_\alpha$ related to $T_\alpha$ in \eqref{eq_torsioncason}, it turns out to be totally skew-symmetric if and only if $q\in\RR$ and $t=\varepsilon q-\frac{n+1}{n}-2\varepsilon $.
In this case,
\begin{equation}\label{eq_skewcason}
(\alpha-\alpha_{g_\varepsilon  })(X,Y)=(1-q)\,\big(\varepsilon (bz-aw),\mathbf{i}\,  \mathrm{Im}( \,\overline{z}^tw)\big)
\end{equation}
and
\begin{equation*}\label{eq_3formacason}
\omega_\alpha(X,Y,Z)=
2\varepsilon (1-q)\Re(au^t\bar w+bz^t\bar u+cw^t\bar z).
\end{equation*}
Following Eq.~\eqref{cur}, we tediously compute the curvature tensor of the map $\alpha$ as in \eqref{eq_skewcason} to get
\begin{equation}\label{eq_curvaturecason}
\begin{array}{ll}\vspace{2pt}
R_\alpha(X,Y,Z)=
 & \big( \frac{\varepsilon q^2}2(z( \bar w^t u- w^t \bar u)+w(z^t\bar u-\bar z^tu))+z(\bar w^t u)-w(\bar z^tu)  \\
 \vspace{3pt}
  &+(-\varepsilon q+2\varepsilon +1)u(\bar w^tz-\bar z^t w)+
   \varepsilon ^2 (q^2-2q)  c(bz-aw), \\
 &\qquad\qquad\  \quad-\frac12\varepsilon  (q^2-2q) ((\bar z^tu+z^t\bar u)b-(\bar w^tu+w^t\bar u)a) \big).
\end{array}
\end{equation}
This allows us to compute the corresponding Ricci tensor
\begin{equation}\label{eq_Riccicason}
\Ric_\alpha(X,Y)=2(\varepsilon(q^2-2q+2)+n+1)\Re(z^t\bar w)   +2n\varepsilon^2(q^2-2q)ab,
\end{equation}
which is always symmetric. Thus, the Ricci tensor  is a scalar multiple of the metric $g_\varepsilon$ if and only if 
$$
\varepsilon(q^2-2q+2)+n+1=n\varepsilon(q^2-2q),
$$
or equivalently if
$$
(q-1)^2=\frac{2\varepsilon+n+1}{\varepsilon(n-1)}+1=\frac{(\varepsilon+1)(n+1)}{\varepsilon(n-1)}.
$$
This gives two solutions whenever $\varepsilon(\varepsilon+1)>0$
 and one solution for $\varepsilon=-1$. This proves   item   $n>3$ in Theorem~\ref{th_main}, simply by taking into account Eq.~\eqref{eq_tensoresenm}, which implies
$$
\Phi(X,Y)\,\xi+\varepsilon(\eta(X)\psi(Y)-\eta(Y)\psi(X))=\big(\varepsilon (bz-aw),\mathbf{i}\,  \mathrm{Im}( \,\overline{z}^tw)\big),
$$
at $o\in\esfera$ under the identification $T_o\esfera\cong\mm\cong\CC^n\oplus\RR\mathbf{i}$.

%%%%%%%%%%%%%%%%%%%%%%%%%%%%%%%%%%%%%%%%%%%%%%%%%%%%%%%%%%%
%%%%%%%%%%%%%%%%%%%%%%%%%%%%%%%%%%%%%%%%%%%%%%%%%%%%%%%%%%%
 \subsection{Case $\esf^7\cong\frac{\SU(4)}{\SU(3)}$}
 
 Again recall from \cite[Eq.~(31)]{DraperGarvinPalomo} that an $\RR$-bilinear map $\beta\colon\mm\times\mm\to\mm$ is $\hh$-invariant if and only if there exist $q_{1},q_{2},q_{3},q_4\in \mathbb{C}$ and $t\in \mathbb{R}$ such that, for any $X,Y\in\mm$,
\begin{equation*}\label{eq_todoscason3}
\beta(X,Y)=\left( q_{1}bz+q_{2}aw+q_4\,\bar z\times\bar w,\, \mathbf{i}\, \left(tab+\mathrm{Im}(q_{3}\overline{z}^tw)\right)\right).
\end{equation*}

The conditions to assure  $\beta\in\mathfrak{so}(\mathfrak{m},g_\varepsilon)$ are $t=\mathrm{Im}(q_{2})=0$ and $q_{1}=\varepsilon\overline{q_{3}}$ as in the previous case, so that, by renaming, there are $q,p\in \mathbb{C}$ and $t\in \mathbb{R}$  such that 
\begin{equation*}\label{eq_metriccason3}
\beta(X,Y)=\left( -\varepsilon q\,bz+t\,aw+p\,\bar z\times\bar w,\, -\mathbf{i}\,  \mathrm{Im}(\bar q\,\overline{z}^tw)\right).
\end{equation*}

The torsion tensor $T{_{\beta}}$ is related to that one in Eq.~\eqref{eq_torsioncason} by
\begin{equation*}\label{eq_torsioncason3}
(T{_{\beta}}-T_\alpha)(X,Y)=(2p\,\bar z\times\bar w,0),
\end{equation*}
which vanishes for $q=1,\, t=-\varepsilon -\frac{n+1}{n},\,p=0$, and the Levi-Civita connection is given again by Eq.~\eqref{eq_LCcason}.

The conditions for having skew-torsion  are still $q\in\RR$ and $t=\varepsilon q-\frac{n+1}{n}-2\varepsilon $, for arbitrary $p$.
In this case,
\begin{equation*}\label{eq_skewcason3}
(\beta-\alpha_{g_\varepsilon  })(X,Y)=s\,\big(\varepsilon (bz-aw),\mathbf{i}\,  \mathrm{Im}( \,\overline{z}^tw)\big)+(p\,\bar z\times\bar w,0),
\end{equation*}
for $s=1-q\in\RR$ and $p\in\CC$,
and the 3-form becomes
\begin{equation*}\label{eq_3formacason3}
\omega_\beta(X,Y,Z)=
2\varepsilon s\,\Re(au^t\bar w+bz^t\bar u+cw^t\bar z)+2\,\Re(\bar p\det(z,w,u)).
\end{equation*}
The curvature tensor is related to that one in Eq.~\eqref{eq_curvaturecason} by
\begin{equation}\label{eq_curvaturecason3}
\begin{array}{ll}
(R_\beta-R_\alpha)(X,Y,Z)=
 & \big(
            (2\varepsilon q-4\varepsilon -4)p\big(a(\bar w\times\bar u)-b(\bar z\times\bar u) \big)+2\varepsilon qpc\,(\bar z\times\bar w) \\
          &  \quad            +p \bar p \,(\bar z\times (w\times u)-\bar w\times (z\times u)),\,2q \mathbf{i}\,\mathrm{Im}(\bar p\det(z,w,u))\big);
 % \\
 %&2q \mathbf{i}\,\mathrm{Im}(\bar p\det(z,w,u))\big).
\end{array}
\end{equation}
and the Ricci tensor  
\begin{equation*} 
(\Ric_\beta-\Ric_\alpha)(X,Y)=-4p\bar p\,\Re(z^t\bar w),
\end{equation*}
 so that Eq.~\eqref{eq_Riccicason} gives
\begin{equation*}\label{eq_Riccicason3}
\Ric_\beta(X,Y)=2(\varepsilon(q^2-2q+2)+4-2p\bar p)\Re(z^t\bar w)   +6\varepsilon^2(q^2-2q)ab,
\end{equation*}
which is always symmetric.
Thus, the Ricci tensor is a scalar multiple of  the metric $g_\varepsilon$ if and only if 
$$
\varepsilon(s^2+1)+4-2p\bar p=3\varepsilon(s^2-1),
$$
or equivalently
$
\varepsilon s^2+p\bar p=2(\varepsilon+1),
$
which gives item $n=3$ in Theorem~\ref{th_main} 
for the choice $ s_1=-\Re(p)$ and $s_2=\mathrm{Im}(p)$,     
by taking into account Eq.~\eqref{eq_tensoresenm3}.
%and  
%$$
%%\begin{array}{l}
% \Re(\det(z,w,u))=g_\varepsilon((\bar z\times\bar w,0) ,Z ),\qquad
%%\\
%  \mathrm{Im}(\det(z,w,u))=g_\varepsilon((\mathbf{i}\bar z\times\bar w,0) ,Z ).
%  % \end{array}
%   $$

 %%%%%%%%%%%%%%%%%%%%%%%%%%%%%%%%%%%%%%%%%%%%%%%%%%%%%%%%%%%
 %%%%%%%%%%%%%%%%%%%%%%%%%%%%%%%%%%%%%%%%%%%%%%%%%%%%%%%%%%%
  \subsection{Case $\esf^5\cong\frac{\SU(3)}{\SU(2)}$}
  
%  Let $\theta\colon\CC^2\to\CC^2$ the $\hh=\mathfrak{su}(2)$-homomorphism given by
%  $
%  \theta(z)=\left(\begin{array}{l}-\bar z_2\\ \bar z_1\end{array}\right).
%$
%For brevity, we denote $X=(z,a),Y=(w,b), Z=(u,c)$ arbitrary elements in $ \mathfrak{m}$.
%Recall from \cite{DraperGarvinPalomo} that an 
Now an arbitrary
$\hh$-invariant bilinear map $\gamma\colon\mm\times\mm\to\mm$ is given by
\begin{equation*}\label{eq_todoscason2}
\gamma(X,Y)=\left( q_{1}bz+q_{2}aw+p_{1}b\theta(z)+p_{2}a\theta(w),\, \mathbf{i}\,  tab+\mathbf{i}\,\mathrm{Im}(q_{3}\overline{z}^tw+p_{3}\overline{\theta(z)}^tw)\right) 
\end{equation*}
%for all $X,Y\in\mm$, 
for some $q_{1},q_{2},q_{3},p_{1},p_{2},p_{3}\in \mathbb{C}$ and $t\in \mathbb{R}$ \cite[\S6.1]{DraperGarvinPalomo},  where the map $\theta\colon\CC^2\to\CC^2$ is defined after Eq.~\eqref{eq_tensoresenm2}.

Those ones compatible with $g_\varepsilon$ form a 7-dimensional vector space, given by the relations
$t=\mathrm{Im}(q_{2})=0,\,q_{1}=\varepsilon\overline{q_{3}},\,p_{1}=\varepsilon\overline{p_{3}}$, so that, by renaming, there are $q,p,p_2\in \mathbb{C}$ and $t\in \mathbb{R}$  such that 
\begin{equation*}\label{eq_metriccason2}
\gamma(X,Y)=\left( -\varepsilon q\,bz+t\,aw-\varepsilon p\,b\theta(z)+p_{2}a\theta(w),\, -\mathbf{i}\,  \mathrm{Im}(\bar q\,\overline{z}^tw+\bar p\,\overline{\theta(z)}^tw)\right).
\end{equation*}
 The torsion tensor $T_\gamma$ is given by
\begin{equation*}\label{eq_torsioncason2}
(T_\gamma-T{_{\alpha}})(X,Y)=\left( (-\varepsilon p-p_2 )(b\theta(z)-a\theta(w)), -2\mathbf{i}\,  \mathrm{Im}( \bar p\,\overline{\theta(z)}^tw) \right),
\end{equation*}
which vanishes for $q=1, t=-\varepsilon -\frac{3}{2},p=p_2=0$. This means that the Levi-Civita connection corresponds again with Eq.~\eqref{eq_LCcason}. Now $\omega_\gamma$ is a 3-form (equivalently, $\gamma-\alpha_{g_\varepsilon}$ is a skew-symmetric map) if and only if $q\in\RR, t=\varepsilon q-\frac{3}{2}-2\varepsilon $ and $p_2=\varepsilon p$.
Hence, an arbitrary $\hh$-invariant bilinear map $\gamma$, metric with skew-torsion, is given by
\begin{equation*}\label{eq_skewcason2}
(\gamma-\alpha_{g_\varepsilon  })(X,Y)= 
%\,\Big(\varepsilon \big(s(bz-aw)-p(b\theta(z)-a\theta(w))\big),\mathbf{i}\,  \mathrm{Im}( \,s\overline{z}^tw-\bar p\,\overline{\theta(z)}^tw)\Big) 
%\end{equation}
%$$
s\,\Big(\varepsilon (bz-aw),\mathbf{i}\,  \mathrm{Im}( \,\overline{z}^tw)\Big)-\Big(\varepsilon p(b\theta(z)-a\theta(w)),\mathbf{i}\,  \mathrm{Im}(  \bar p\,\overline{\theta(z)}^tw)\Big)
\end{equation*}
%$$
for $s=1-q\in\RR$ and $p\in\CC$. For computing the symmetrized Ricci tensor, we can save the computation of the curvature tensor\footnote{In fact, a tedious computation of the curvature shows that in this case, Ricci tensor is not necessarily symmetric.
Alternatively, it is possible to check that $\mathrm{div}(T_\alpha)\ne0$, as done in \cite[Remark~6.10]{DraperGarvinPalomo}.}    by using the formula in \cite[Appendix]{surveyagricola}: 
\begin{equation}\label{formulicas}
\mathrm{Sym}(\mathrm{Ric}_{\alpha})=\mathrm{Ric}_{\alpha_{g_\varepsilon}}-\frac{1}{4}S, 
%\mathrm{Sym}(\mathrm{Ric}^{\nabla})=\mathrm{Ric}^{g}-\frac{1}{4}S, 
%\qquad\qquad s^{\nabla}=s^{g}-\frac{3}{2}\|T^{\nabla}\|^2,
\end{equation}
where the tensor $S\in \mathcal{T}^{0,2}(M)$  at $p\in M$  is defined by
$$
S(X,Y)_p=\sum_{j=1}^{5}g_\varepsilon(T^{\nabla}(e_{j},X_p),T^{\nabla}(e_{j},Y_p))\,g_\varepsilon(e_j,e_j),
$$
for $\{e_{1},...,e_{5}\}$   any orthonormal basis of $T_pM$. Again, by invariance, we work at $p=o$, chosing, for instance, the following orthonormal basis of $\mm$,
$$
e_1=((1,0),0),\,e_2=((\mathbf{i},0),0),\,e_3=((0,1),0),\,e_4=((0,\mathbf{i}),0),\,e_5=\frac1{\sqrt{|\varepsilon|}}((0,0),\mathbf{i}).
$$
Now we check that
$$
S(X,Y)=
%(4\varepsilon^2-4\varepsilon) (s^2+p\bar p)   \Re(z^t\bar w)-16\varepsilon^2    (s^2+p\bar p) ab.MAL
-8\varepsilon (s^2+p\bar p)   \Re(z^t\bar w)-16\varepsilon^2    (s^2+p\bar p) ab.
$$
By extracting from Eq.~\eqref{eq_Riccicason} the Ricci tensor of the Levi-Civita connection (making $q=1, n=2$),
\begin{equation*} 
\begin{array}{ll}
\Ric_{\alpha_{g_\varepsilon}}(X,Y)&=2(\varepsilon +n+1)\Re(z^t\bar w)   -2n\varepsilon^2 ab\\
&=2(\varepsilon +3)\Re(z^t\bar w)   -4\varepsilon^2 ab,
\end{array}
\end{equation*}
we can substitute in Eq.~\eqref{formulicas} to get
\begin{equation*}
\mathrm{Sym}(\mathrm{Ric}_{\gamma})=\big(2\varepsilon+6+2\varepsilon  (s^2+p\bar p) \big)\,  \Re(z^t\bar w)+    (s^2+p\bar p-1) \,4\varepsilon^2 ab. 
%\qquad\qquad s^{\nabla}=s^{g}-\frac{3}{2}\|T^{\nabla}\|^2,
\end{equation*}
Thus, the symmetric part  of the Ricci tensor  is a scalar multiple of    the metric $g_\varepsilon$ if and only if 
$$
2\varepsilon+6+2\varepsilon  (s^2+p\bar p)=4\varepsilon(s^2+p\bar p-1), \,
$$
or equivalently
$
 s^2+p\bar p=3\frac{\varepsilon+1}\varepsilon.
$
This equation has solutions whenever $({\varepsilon+1})\varepsilon>0$, which gives item $n=2$ in Theorem~\ref{th_main} by taking $ s_3=-\Re(p)$ and $s_4=\mathrm{Im}(p)$,     
taking into account 
 Eqs.~\eqref{eq_tensoresenm2}, \eqref{eq_tensoresenm} and noting that       
$$
\begin{array}{ll}
 \Phi(\hat\psi(X),Y)=(0, \mathbf{i}\, \mathrm{Im}(    \overline{\theta(z)}^tw)),&
 \eta(X)\psi(\hat\psi(Y))-\eta(Y)\psi(\hat\psi(X)) =( b\theta(z)-a\theta(w),0),\\
 \Phi(\psi(\hat\psi(X)),Y)=(0, -\mathbf{i}\, \mathrm{Re}(    \overline{\theta(z)}^tw)) ,
  &
  \eta(X)\hat\psi(Y)-\eta(Y)\hat\psi(X)=( -\mathbf{i}\,  (b\theta(z)-a\theta(w)),0).\\
 \end{array}
$$
%at $o\in\esf^5$ under the identification $T_o\esf^5\cong\mm\cong\CC^2\oplus\RR\mathbf{i}$.

We also observe   $\Phi(\psi(\hat\psi(X)),Y)=g_\varepsilon( \hat\psi(X),Y)$ (a not closed invariant 2-form).

  %%%%%%%%%%%%%%%%%%%%%%%%%%%%%%%%%%%%%%%%%%%%%%%%%%%%%%%%%%%
  %%%%%%%%%%%%%%%%%%%%%%%%%%%%%%%%%%%%%%%%%%%%%%%%%%%%%%%%%%%
   \subsection{Case $\esf^3\cong \SU(2) $}

We know from Eq.~\eqref{eq_skewcason} that the map $\alpha$ defined by
\begin{equation}\label{eq_skewcason1}
(\alpha-\alpha_{g_\varepsilon  })(X,Y)=(1-q)\,\big(\varepsilon (bz-aw),\mathbf{i}\,  \mathrm{Im}( \,\overline{z}^tw)\big),
\end{equation}
provides a metric with skew-torsion connection for any $q\in\RR$,
where recall that the Levi-Civita connection is given through Nomizu's theorem by
 \begin{equation*}\label{eq_LCcason1}
\alpha_{g_\varepsilon  }(X,Y)=\left( -\varepsilon  \,bz-(\varepsilon +2)\,aw,\,
 -\mathbf{i}\,  \mathrm{Im}( \,\overline{z}^tw)\right).
\end{equation*}
Now, taking into account that, in this case,  $\wedge^3\mm\cong\RR$, and hence $\hom(\wedge^3\mm,\RR)$ is one-dimensional ($\hh=0$), the family in Eq.~\eqref{eq_skewcason1} exhausts the set of maps related to metric with skew-torsion invariant affine connections.

Then Eq.~\eqref{eq_Riccicason} provides the Ricci tensor
\begin{equation*}\label{eq_Riccicason1}
\Ric_\alpha(X,Y)=2(\varepsilon(q^2-2q+2)+2)\Re(z^t\bar w)   +2\varepsilon^2(q^2-2q)ab,
\end{equation*}
which is a scalar multiple of the metric $g_\varepsilon$ if and only if 
$
\varepsilon(q^2-2q+2)+2=\varepsilon(q^2-2q),
$
that is, if and only if $\varepsilon=-1$.\smallskip

%%%%%%%%%%%%%%%%%%%%%%%%%%%%%%%%%%%%%%%%%%%%%%%%%

%\textcolor{blue}{\section{Some comments and conclusions QUIZA}}
%\begin{table}\label{tabla}
%\centering
%\begin{tabular}{ c||c|c|c|c||}
%  & $n\ge3$ & $n=3$ & $n=2$ & $n=1$ \\
%\hline\hline
%   $\varepsilon\in(-\infty,-1)$ & $ 2$ pt. &hyperbolloid   2-sh& ellipsoid & $\emptyset$ \\
%\hline
%$\varepsilon=-1$ & $ 1$ pt.&cone& 1 pt.& line \\
%\hline
%$\varepsilon\in( -1,0)$ & $\emptyset$ &hyperbolloid   1-sh&  $\emptyset$ &$\emptyset$\\
%\hline
%\hline
%$\varepsilon>0$ & $2$ pt.&ellipsoid& ellipsoid &$\emptyset$\\
%\hline\hline
%\end{tabular}\vspace{3pt}
%\caption{Varieties parametrizing Einstein connections with skew-torsion}\label{tablatita}
%\end{table} 

\medskip

%\newpage

We have finished the proof of  Theorem~\ref{th_main}. The results are compiled in the following table.\bigskip

%\begin{table}%\label{tabla}
%\centering
\begin{center}
\begin{tabularx}{14cm}{| X c  | c | c| c | c |}
%\begin{tabularx}{15cm}{| c  l | c | X |}
&& \ $n\ge4$ \ & \ $n=3$\  &\  $n=2$ \ &\  $n=1$\  \\
\hline \hline
 &$\varepsilon\in(-\infty,-1)$ & $ 2$ pt. &\ hyperboloid   2-sheets\ & \ ellipsoid\  & $\emptyset$
  \\  %\hdashline
\cdashline{2-6} \vrule width 0pt height 14pt
 Riemann& $\varepsilon=-1$ & $ 1$ pt.&cone& 1 pt.& line \\
 %\hline 
 \cdashline{2-6} \vrule width 0pt height 14pt
 & $\varepsilon\in( -1,0)$ & $\emptyset$ &hyperboloid   1-sheet&  $\emptyset$ &$\emptyset$\\
 %\hline 
%
%\end{tabularx}
%
%\begin{tabularx}{\textwidth}{ |c|l|c|X|}
\hline\hline
Lorentz &$\varepsilon>0$ & $2$ pt.&ellipsoid& ellipsoid &$\emptyset$\\
\hline   
\end{tabularx}
\vspace{3pt}

%\caption
\textsc{Table 1}. {\small  Varieties parametrizing invariant Einstein  with skew-torsion connections.}\label{tablatita}
\end{center}\vspace{3pt}
%\end{table} 
 
\begin{remark}\rm{
As a consequence of Eq.~\eqref{eq_Einsteincondition},
% the previous formulae for $\mathrm{Sym}(\mathrm{Ric}_{\alpha})$, 
and with the notations used in Theorem~\ref{th_main}, the scalar curvature of an Einstein  with skew-torsion invariant connection becomes
%$s^\nabla=\frac{2n}{2n+1}\varepsilon(s^2-1)$ if $n\ne2$, while $ s^\nabla=\frac{4}{5}\varepsilon(s^2+s_3^2+s_4^2-1)$ if $n=2$.
$$
 s^\nabla=\left\{
 \begin{array}{ll}\vspace{3pt} {2n}{(2n+1)}\,\varepsilon(s^2-1)&\textrm{ if $n\ne2,$}\\
    20\,\varepsilon(s^2+s_3^2+s_4^2-1) &\textrm{ if $n=2$.}\end{array}\right.
$$
Hence, by taking into consideration the relationship among the parameters,
$$
 s^\nabla=\left\{
 \begin{array}{ll}\vspace{2pt} {2n}{(2n+1)}\,\frac{2\varepsilon+n+1}{n-1}&\textrm{ if either $n=2$ or $n\ge4$,}\\
 \vspace{2pt}
    6\,(1-s^2) &\textrm{ if $n=1$,}\\
    42\,(\varepsilon+2-s_1^2-s_2^2) &\textrm{ if $n=3$.}\end{array}\right.
$$

Thus, there are Ricci-flat invariant Einstein  with skew-torsion connections    for any odd-dimensional sphere for a suitable Berger metric, which are precisely,
\begin{itemize}
\item For $n=1$, the round metric ($\varepsilon=-1$) and  $s=\pm1$;
\item For $n=2$, the metric associated to $\varepsilon=-3/2$ and $s^2+s_3^2+s_4^2=1$;
\item For $n=3$, any metric $g_\varepsilon$ with $0\ne\varepsilon\ge-2$, and $s=\pm1,s_1^2+s_2^2=\varepsilon+2$;
\item For $n\ge4$, the metric associated to  $\varepsilon=-\frac{n+1}2$ and $s=\pm1$.
\end{itemize}
In particular, only in   $\esf^7$ there are Lorentzian Berger spheres satisfying such properties. In fact, in such a case every Lorentzian sphere does.

Also, we can conclude from Eq.~\eqref{eq_curvaturecason} that, for $n\ge4$, there do not exist any flat invariant affine connection    with skew-torsion (well-known from the work by Cartan and Schouten \cite{CartanSchouten}), while for $n=3$, Eq.~\eqref{eq_curvaturecason3} implies that there are flat $\SU(4)$-invariant affine connections   with skew-torsion in $\esf^7$ only for the round metric, and, in that case, the flat invariant connections correspond to $s=1$ and arbitrary $s_1,s_2$ such that $s_1^2+s_2^2=1$.}

\end{remark}

%\vspace{1truecm}

\end{document}